\documentclass{article}
\usepackage[applemac]{inputenc}
\usepackage[tbtags]{amsmath}
\usepackage{hyperref}
\usepackage{xcolor}
\usepackage[pdftex]{graphicx} 
\usepackage{tikz}
\DeclareRobustCommand\apr{\mathrel{\mathchar"1270\mkern2mu}}

\newcommand\rhoi[1]{\rho\raise.5ex\hbox{${}_{[#1]}$}}

\newcommand\xk{x_k}

\newcommand\pq[2]{\raise.25ex\hbox{\footnotesize${#1}\over{#2}$}%
\hskip-.35ex\null}
\newcommand\pqbis[2]{\leavevmode\raise.5ex\hbox{\scriptsize$#1$}\hskip-.45ex\raise.25ex\hbox{\scriptsize$/$}\hskip-.4ex\lower.2ex\hbox{\scriptsize$#2$}}
\newcommand\onehalf{\leavevmode\raise.5ex\hbox{\scriptsize$1$}\hskip-.25ex/\hskip-.25ex\lower.2ex\hbox{\scriptsize$2$}}

\newcount\npar
\newcommand\nouparc{\advance\npar by1\medskip\noindent\hbox to\parindent{\textbf{\the\npar.}}}

\begin{document}

\noindent
\textbf{Phragmén’s sequential method with a variance criterion}

\smallskip
\noindent
Xavier Mora, Dept. Matemàtiques, Univ. Autònoma de Barcelona,\hfil\break
7th September 2016

\bigskip
This short note briefly explores (a variant of) a method for electing a parliament by means of approval voting. Phragmén's standard sequential method \cite{janson:2012,janson:2016,mora:2015} is a method of this kind that generalizes D'Hondt's rule: every step proceeds so as to minimize the maximum representation. What about similarly generalizing Sainte-Laguë's rule, i.\,e.\ minimizing the variance at every step?

\nouparc
We'll follow the notation of \cite{mora:2015}. Assume that $n$ seats have already been assigned, and that each elector of type $k$ is represented in the amount $r_k$. We~consider allocating a new seat to candidate (or party) $i$ and distributing it between the electors that approved $i$. Let $\xk$ be the amount of this new seat that is given to each elector of type $k$ (we are omitting the reference to $i$). The $\xk$ are subject to the following constraints:
\begin{alignat}{2}
\label{eq:positivenessconstraint}
\xk &\ge 0,\qquad&&\text{for all $k$}
\\[2pt]
\label{eq:zeroconstraint}
\xk &= 0,\qquad&&\text{for $k\not\apr i$}
\\[2pt]
\label{eq:oneconstraint}
\sum_{k\apr i} u_k \xk &= 1,\qquad&&
\end{alignat}

The variance (multiplied by $w$) of the new distribution is
\begin{multline}
V = \sum_k u_k \left(r_k+\xk-\frac{n+1}{w}\right)^2 = \sum_k u_k \left(r_k+\xk\right)^2 - \frac{(n+1)^2}{w}\\
= \sum_k u_k r_k^2 + \sum_{k\apr i} u_k \left(2r_k\xk+\xk^2\right) - \frac{(n+1)^2}{w},
\qquad\null
\end{multline}
where we have used \eqref{eq:zeroconstraint}. So it suffices to minimize the value of
\begin{equation}
\label{eq:phi}
\phi_i = \frac12\, \sum_{k\apr i} u_k \left(2r_k\xk+\xk^2\right),
\end{equation}
We will now minimize this value for each candidate $i$ and then will minimize the result with respect to $i$.

Introducing a Lagrange multiplier $\rho$ to deal with the constraint \eqref{eq:oneconstraint}, we have to solve the following system of (linear) equations, where the unknonws are $\xk$ for $k\apr i$ and $\rho$:
\begin{equation}
\label{eq:derivative}
\partial\phi_i/\partial \xk - \rho u_k = u_k (r_k + \xk) - \rho u_k = 0.
\qquad\text{for each $k\apr i,$}
\end{equation}
together with equation \eqref{eq:oneconstraint}.
Summing \eqref{eq:derivative} for all $k\apr i$ and using \eqref{eq:oneconstraint} one gets
\begin{equation}
\bigg(\sum_{k\apr i} u_k r_k\bigg) + 1 - \rho w_i = 0,
\end{equation}
where we are using the notation $w_i = \sum_{k\apr i} u_k,$ as in \cite{mora:2015}. So we get
\begin{equation}
\label{eq:lambdai}
\displaystyle
\rho = \frac{\big(\sum_{k\apr i} u_k r_k\big) + 1}{w_i} =: \rho_i,
\end{equation}
which introduced in \eqref{eq:derivative} results in the new representations being
\begin{equation}
\label{eq:samevalue}
\displaystyle
r_k + \xk = \rho_i = \frac{\big(\sum_{k\apr i} u_k r_k\big) + 1}{w_i},\qquad\text{for each $k\apr i,$}
\end{equation}
where the right-hand side does not depend on $k\apr i$. Notice that this value is exactly the one that appears in equation (29) of \cite{mora:2015}.

Now we have to minimize $\phi_i$ with respect to $i$. By introducing \eqref{eq:samevalue} in \eqref{eq:phi}, we get 
\begin{equation}
\label{eq:phii}
2\,\phi_i = \sum_{k\apr i} u_k (\rho_i-r_k)(\rho_i+r_k)
= \sum_{k\apr i} u_k (\rho_i^2-r_k^2)
= w_i\rho_i^2 -  \sum_{k\apr i} u_k r_k^2,
\end{equation}
where $\rho_i$ is the value given by \eqref{eq:lambdai}.

\smallskip
So, it is a matter of \emph{choosing the $i$ that minimizes \eqref{eq:phii} and distributing this new seat according to \eqref{eq:samevalue}}. This is valid under the assumption that the obtained $\xk$ satisfy \eqref{eq:positivenessconstraint}. Unfortunately, it is not always so, which issue will be considered in a while.

\nouparc
In the case of closed party lists this algorithm reduces indeed to Sainte-Laguë's rule. In fact, every elector who approves party $i$ has a representation $r_k = n_i/w_i$. And their total representation is $\sum_{k\apr i}u_kr_k = n_i$ (since $\sum_{k\apr i}u_k = w_i$). By plugging this into \eqref{eq:lambdai}, we get $\rho_i = (n_i+1)/w_i$. And by plugging all this into \eqref{eq:phii}, we get
\begin{equation}
\label{eq:phiiSL}
2\,\phi_i = \frac{(n_i+1)^2}{w_i} - \frac{n_i^2}{w_i} = \frac{2 n_i+1}{w_i},
\end{equation}
whose minimization gives Sainte-Laguë's algorithm.

\nouparc
As we have already said, sometimes the above obtained distribution does not satisfy the positiveness constraint \eqref{eq:positivenessconstraint} (in contrast to the Phragmén's standard sequential method, where positivity is guaranteed by \cite[Prop.\,7.4]{mora:2015}). An example of such a failure is provided by the following profile
\begin{equation}
9\,:\ \textsf{a1,\,a2},\quad 1\,:\ \textsf{a1,\,a2,\,b},\quad 3\,:\ \textsf{b,\,c}.
\end{equation}
In this case, for $n = 3$ the preceding computations lead to the following assignment:

\def\bstrut{\rule[-.85ex]{0ex}{3.25ex}}
\begin{center}
\begin{tabular}{|c|c|c|c|c|c|}
\hline\bstrut
Seat \# & $i$ & 9\,:\,\textsf{a1,\,a2} & 1\,:\,\textsf{a1,\,a2,\,b} & 3\,:\,\textsf{b,\,c}\\\hline\bstrut
1 & \textsf{a1} & 0.1000  & 0.1000 & \textcolor{gray}{0.0000}\\
\bstrut
2 & \textsf{b} & \textcolor{gray}{0.0000}  & 0.1750 & 0.2750\\
\bstrut
3 & \textsf{a2} & 0.1175  & \llap{$-$}0,0575 & \textcolor{gray}{0.0000}\\
\hline
\end{tabular}
\end{center}

Such cases can be corrected by explicitly imposing a constraint $x_k = 0$ whenever the above computations lead to a negative value of $x_k$, as if the electors $k$ hadn't approved candidate $i$. In the preceding example, this results in the following assignment

\begin{center}
\begin{tabular}{|c|c|c|c|c|}
\hline\bstrut
Seat \# & $i$ & 9\,:\,\textsf{a1,\,a2} & 1\,:\,\textsf{a1,\,a2,\,b} & 3\,:\,\textsf{b,\,c}\\
\hline\bstrut
1 & \textsf{a1} & 0.1000  & 0.1000 & \textcolor{gray}{0.0000}\\
\bstrut
2 & \textsf{b} & \textcolor{gray}{0.0000}  & 0.1750 & 0.2750\\
\bstrut
3 & \textsf{a2} & 0.1111 & \textcolor{gray}{0.0000} & \textcolor{gray}{0.0000}\\
\hline
\end{tabular}
\end{center}

Anyway, it remains to be checked whether such a correction will always give the exact minimizer of the original problem.

\nouparc
Like the standard sequential Phragmén's method, the party version of the present variation also lacks monotonicity. The following example is similar to (48) of~\cite[\S7.5]{mora:2015}:
\begin{equation}
4\,:\ \textsf{A},\quad 1\,:\ \textsf{B},\quad 3\,:\ \textsf{C},\quad 9\,:\ \textsf{A,\,B},\quad 3\,:\ \textsf{B,\,C}.
\end{equation}
For $n = 3$ this profile gives 2~seats to \textsf{A}, whereas adding one vote that approves only~\textsf{A}, i.\,e. changing the coefficients to $(5,1,3,9,3)$, results in this party getting only 1~seat.

\nouparc
On the other hand, for two parties, the limit $n\rightarrow\infty$ shows a Cantorian staircase phenomenon similar to that of \cite[\S7.7]{mora:2015}. For instance, the next figure shows the dependence of $f_{1200}$ on $\alpha$ for $\zeta = 0.376$.

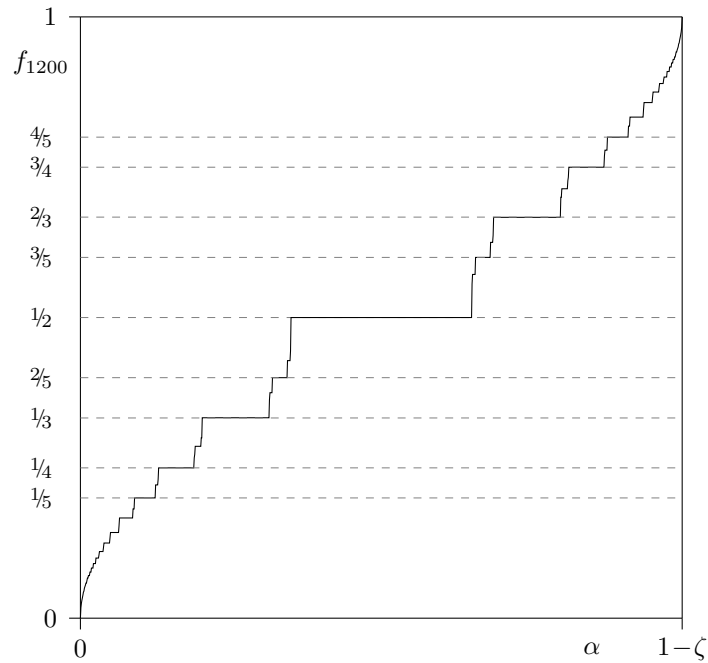
\begin{center}
\begin{tikzpicture}[x=80mm,y=80mm,>=stealth]
\draw (0,0) -- (1,0);
\node at (0.85,-4mm) {$\alpha$};
\draw (0,0) -- (0,1);
\node at (-5.25mm,0.925) {$f_{1200}$};
\draw (1,0) -- (1,1);
\draw (0,1) -- (1,1);
\draw (-1.5mm,0) -- (0,0);
\node at (-4mm,0) {$0$};
\draw (-1.5mm,1) -- (0,1);
\node at (-4mm,1) {$1$};
\draw (0,-1.5mm) -- (0,0);
\node at (0,-4mm) {$0$};
\draw (1,-1.5mm) -- (1,0);
\node at (1,-4mm) {$1\!-\!\zeta$};
\draw[help lines, dashed] (0,0.5) -- (1,0.5);
\node at (-5.25mm,0.5) {$\pqbis12$};
\draw plot file {replansSL1200.txt};
\draw[help lines, dashed] (0,0.333) -- (1,0.333);
\node at (-5.25mm,0.333) {$\pqbis13$};
\draw[help lines, dashed] (0,0.667) -- (1,0.667);
\node at (-5.25mm,0.667) {$\pqbis23$};
\draw[help lines, dashed] (0,0.25) -- (1,0.25);
\node at (-5.25mm,0.25) {$\pqbis14$};
\draw[help lines, dashed] (0,0.75) -- (1,0.75);
\node at (-5.25mm,0.75) {$\pqbis34$};
\draw[help lines, dashed] (0,0.2) -- (1,0.2);
\node at (-5.25mm,0.2) {$\pqbis15$};
\draw[help lines, dashed] (0,0.4) -- (1,0.4);
\node at (-5.25mm,0.4) {$\pqbis25$};
\draw[help lines, dashed] (0,0.6) -- (1,0.6);
\node at (-5.25mm,0.6) {$\pqbis35$};
\draw[help lines, dashed] (0,0.8) -- (1,0.8);
\node at (-5.25mm,0.8) {$\pqbis45$};
\end{tikzpicture}%
\end{center}

\nouparc
Summing up, except for the fact of its being indeed a generalization of Sainte-Laguë's rule, this method doesn't seem especially interesting in comparison with Phragmén's standard sequential method. The rather undesirable phenomena of the preceding paragraphs are still present. And the fact that positivity is not immediately satisfied is quite annoying.


\begin{thebibliography}{99\kern2pt}
\bgroup\small 

\newcommand\ensep{\hskip.55em plus.25em}
\newcommand\tita[1]{{\upshape #1}}
\newcommand\titl[1]{{\itshape #1}}
\newcommand\tits[1]{{\itshape #1}}
\newcommand\vvv{\unskip, }
\newcommand\nnn{\unskip, n.\,}
\newcommand\ppp{\unskip\,: }

\bibitem{janson:2012}
Svante Janson, 2012--2015.\ensep
\titl{Proportionella Valmetoder}.\ensep
\url{http://www2.math.uu.se/~svante/papers/sjV6.pdf}.

\bibitem{janson:2016}
Svante Janson, 2016.\ensep
\tita{Phragmén's and Thiele's election methods}.\ensep
In preparation.

\bibitem{mora:2015}
Xavier Mora, Maria Oliver, 2015.\ensep
\tita{Eleccions mitjançant el vot d'aprovació. El mètode de Phragmén i algunes variants}.\ensep
\tits{Butlletí de la Societat Catalana de Matemàtiques} \vvv30 \ppp57--101.

\egroup\small 
\end{thebibliography}
\end{document}